\documentclass[11pt,a4paper]{article}

\usepackage{amsmath, amssymb, amsthm, amsfonts}
\usepackage{mathrsfs}

\theoremstyle{plain}
\newtheorem{thm}{Theorem}
\newtheorem{cor}{Corollary}
\newtheorem{lem}{Lemma}

\theoremstyle{remark}

\renewcommand{\Re}{{\rm Re\,}}

\renewcommand{\a}{\alpha}

\newcommand{\ba}{\boldsymbol{a}}
\renewcommand{\d}{\displaystyle}
\newcommand{\Psum}{\mathop{{\;\,{\sum}^{\prime}}}}
\newcommand{\Jsum}{\mathop{{\sum_{l=0}^J \sum_{m=0}^J}}}
\newcommand{\Dsum}{\mathop{{\sum \sum}}}

\numberwithin{equation}{section}

\begin{document}

\title{Mean square of the error term in the asymmetric many dimensional divisor problem}
\author{Xiaodong Cao,  Yoshio Tanigawa and Wenguang Zhai}
\date{}

\maketitle

\footnote[0]{2010 Mathematics Subject Classification:11N37}
\footnote[0]{Key words and phrases: Asymmetric many dimensional divisor problem, Mean square of the error term,
Dirichlet series, Functional equation, the Tong-type representation}

\footnote[0]{The first and the third authors are supported by the National Key Basic Research Program of China (Grant No.2013CB834201),
the National Natural Science Foundation of China (Grant No.11171344), the Natural
Science Foundation of Beijing (Grant No.1112010) and the Fundamental Research Funds for the
Central Universities in China (2012Ys01). The second author is supported by Grant-in-Aid for Scientific Research no.24540015.}

\begin{abstract}
Let $\ba=(a_1,a_2,\ldots,a_k)$, where $a_j \ (j=1,\ldots,k)$ are positive integers such that $a_1 \leq a_2 \leq \cdots \leq a_k$.
Let $d(\ba;n)=\sum_{n_1^{a_1}\cdots n_k^{a_k}=n}1$ and $\Delta(\ba;x)$ be the error term of the summatory function of $d(\ba;n)$.
In this paper we show an asymptotic formula of the mean square of $\Delta(\ba;x)$ under a certain condition. Furthermore, in the
cases $k=2$ and 3, we give unconditional asymptotic formulas for these mean squares.
\end{abstract}

\section{Introduction and the statement of results}

Let $k$ be a fixed positive integer and $x \geq 1$. Let $\ba=(a_1,\ldots,a_k)$, where $a_j \ (j=1,\ldots,k)$
are positive integers such that $a_1 \leq \cdots \leq a_k$. Let $d(\ba;n)$ denote the number of representations of
an integer $n$ in the form $n=n_1^{a_1}\cdots n_k^{a_k}$;
\begin{equation} \label{d-def}
d(\ba;n)=\sum_{n_1^{a_1}\cdots n_k^{a_k}=n} 1.
\end{equation}
We define the error term
\begin{equation*}
\Delta(\ba;x):=\Psum_{n \leq x} d(\ba;n)-H(\ba;x),
\end{equation*}
where $H(\ba;x)$ is the main term of the summatory function of $d(\ba;n)$, which is given by the sum of residues of $\prod_{j=1}^k \zeta(a_js)\frac{x^s}{s}$,
and $'$ in the summation symbol means that the last term $d(\ba;x)$ should be halved when $x$ is an integer.
Asymmetric many dimensional divisor problem (or the general divisor problem) is to study the behaviour of $\Delta(\ba;x)$.
See also Ivi\'c \cite{I1} and Kr\"atzel \cite{Kra}, or the survey paper \cite{IKKN}.

When $a_1=a_2=1$, $d(1,1;n)=\sum_{d|n}1, \  \Delta(1,1;x)= \sum_{n \leq  x}d(1,1,;n)-x(\log x+2\gamma-1)$,
($\gamma$ is the Euler  constant),  we have the classical Dirichlet divisor problem.
Dirichlet proved $\Delta(1,1;x)=O(x^{1/2})$ by his famous hyperbola method. The exponent $1/2$ was improved
by many researchers. The latest result is
$$\Delta(x)=O(x^{131/416}(\log x)^{26947/8320})$$
due to Huxley \cite{Huxley}. For the lower bounds, it is known that
$$
\Delta(1,1;x)=\Omega_{+}\left(x^{\frac14}(\log x)^{\frac14}(\log\log x)^{\frac{3+\log 4}{4}}\exp(-c\sqrt{\log\log\log x})\right) \ (c>0)
$$
and
$$
\Delta(1,1;x)=\Omega_{-}\left(x^{\frac14}\exp(c'(\log\log x)^{\frac14}(\log\log\log x)^{-\frac34})\right) \ (c'>0),
$$
which are due to Hafner \cite{Hafner} and Corr\'adi and K\'atai \cite{CK}, respectively. Many corresponding results for asymmetric many dimensional divisor
problem can be found in \cite{I1} and \cite{Kra}.

The mean square estimate is one of the main problems in the theory of divisor problem.
Let $R(T)$ be the error term defined by
$$
\int_1^T \Delta^2(1,1;x)dx=cT^{3/2}+R(T),
$$
where $c$ is a constant given by $c=\frac{1}{6\pi^2}\sum_{n=1}^{\infty}\frac{d(1,1;n)^2}{n^{3/2}}$. Cram\'er \cite{C} first proved that
\begin{equation*} 
R(T)=O(T^{5/4+\varepsilon}).
\end{equation*}
The above estimate of $R(T)$ was improved to
\begin{equation} \label{Tong}
R(T)=O(T\log^5T)
\end{equation}
by Tong \cite{Tong} and recently to $R(T)=O(T\log^3T\log\log T)$ by Lau and Tsang \cite{LT}. Tong's method which derives \eqref{Tong}
is the initial motivation of our previous paper \cite{CTZ}.

Ivi\'c \cite{I2}  studied the upper bound and $\Omega$-result of the mean square of $\Delta(\ba;x)$ for general $k$.
As for the upper bound, he proved that if
$$
\int_1^T \Delta^2(\ba;x) dx \ll T^{1+2\beta_k}   \ \ (\beta_k \geq 0)
$$
then $\beta_k \geq g_k$, where
$$
g_k=\frac{r-1}{2(a_1+\cdots +a_r)}
$$
and $r$ is the largest integer such that
$$
(r-2)a_r \leq a_1+\cdots+a_{r-1} \quad (2 \leq r \leq k)
$$
\cite[(1.5)]{I2}. Furthermore, he showed that the estimate $\int_1^T|\zeta(1/2+it)|^{2k-2}dt \ll T^{1+\varepsilon}$ implies $\beta_k=g_k$.
In particular, $\beta_k=g_k$ holds for $k=2$ and $3$. For the lower bound, he showed that
$$
\int_1^T\Delta^2(\ba;x)dx = \Omega(T^{1+2g_k}\log^A T)
$$
with some constant $A\geq 0$. From these evidence, he conjectured that
\begin{equation}  \label{IvicConj}
\int_1^T \Delta^2(\ba;x) dx =(E_k+o(1))T^{1+2g_k} \log ^A T
\end{equation}
for general $k$  with some constants $E_k>0$ and  $A \geq 0$ \cite[(5.7)]{I2}.

Ivi\'c's conjecture \eqref{IvicConj} was solved by Cao and Zhai \cite{ZC} in the case $k=2$.
More precisely they proved that
\begin{equation} \label{ZhaiCao}
\int_1^T \Delta^2(\ba;x)dx=c(\ba) T^{\frac{1+a_1+a_2}{a_1+a_2}}+O\left(T^{\frac{1+a_1+a_2}{a_1+a_2}-\frac{a_1}{2a_2(a_1+a_2)(a_1+a_2-1)}}\log^{\frac72}T\right),
\end{equation}
where $a_1$ and $a_2$ are integers such that $1 \leq a_1 \leq a_2$, $\ba=(a_1,a_2)$ and $c(\ba)$ is some constant.
Their method is based on the transformation formula of the exponential sum applied to Chowla and Walum type representation
of $\Delta(\ba;x)$ (see also \cite{CTZ1}).
When $a_1=a_2=1$, the error term in \eqref{ZhaiCao} becomes $O(T^{\frac54}\log^{\frac72}T)$. Hence it can be said that
\eqref{ZhaiCao} corresponds to the result of Cram\'er.

In this paper we shall study the mean square estimate of the error term $\Delta(\ba;x)$ more closely by the Tong's method \cite{CTZ, Tong}.
For this purpose, we need an auxiliary divisor function defined by
\begin{equation} \label{dhat-def}
\hat{d}(\ba;n)=\sum_{n_1^{a_1}\cdots n_k^{a_k}=n} n_1^{a_1-1}\cdots n_k^{a_k-1}.
\end{equation}
This function can be  regarded as a dual function of  $d(\ba;n)$.
For notational convenience we put
$$
b(n)=\pi^{2\a-k/2}\hat{d}(\ba;n)  \quad \mbox{and} \quad \mu_n=\pi^{2\a}n,
$$
where
$$
\a=(a_1+\cdots+a_k)/2.
$$

By \eqref{d-def} and \eqref{dhat-def},  it is easy to see that
\begin{equation*}
\varphi(s):= \sum_{n=1}^{\infty}\frac{d(\ba;n)}{n^s}=\prod_{j=1}^k \zeta(a_j s) \qquad  \Re s>1/a_1
\end{equation*}
and
\begin{align}
\psi(s):=& \sum_{n=1}^{\infty}\frac{b(n)}{\mu_n^s} = \pi^{2\a-k/2-2\a s}\sum_{n=1}^{\infty}\frac{\hat{d}(\ba;n)}{n^s} \notag \\
        =& \pi^{2\a-k/2-2\a s}\prod_{j=1}^k \zeta(a_j s-a_j+1)  \qquad \Re s>1.  \label{def-psi}
\end{align}

Let $1/2 \leq \sigma^{\ast}<1 $ be a real number defined by
\begin{equation} \label{def-sigmastar}
\sigma^{\ast}:=\inf\left\{\sigma \ \Big| \  \int_0^T |\psi(\sigma+it)|^2dt \ll T^{1+\varepsilon} \right\}
\end{equation}
for any $\varepsilon>0$. Clearly we have
\begin{equation} \label{sigmastar-lowerbound}
\sigma^{\ast} \geq 1-\frac{1}{2a_k}.
\end{equation}
In this paper we suppose that $\sigma^{\ast}$ satisfies the condition
\begin{equation}  \label{sigmastar-katei}
\sigma^{\ast}<1-\frac{k-1}{4\a}.
\end{equation}
This condition plays an important role in Tong's theory. From \eqref{sigmastar-lowerbound}, we note that \eqref{sigmastar-katei} implies,
as a necessary condition, that
\begin{equation}
(k-2)a_k <a_1+\cdots+a_{k-1}.  \label{a-nokatei}
\end{equation}

We first show an asymptotic formula of the mean square of $\Delta(\ba,x)$.

\begin{thm}
Let $(k-2)a_k <a_1+\cdots+a_{k-1}$ and suppose that $\sigma^{\ast}$ satisfies \eqref{sigmastar-katei}. Then we have
\begin{equation}  \label{maintheorem}
\int_1^T \Delta^2(\ba;x)dx=c(\ba)T^{1+\frac{k-1}{2\a}}+O\left(T^{1+\frac{k-1}{2\a}-\eta(\ba)+\varepsilon}\right),
\end{equation}
where $c(\ba)$ is a certain positive constant and
\begin{equation} \label{eta}
\eta(\ba)=\frac{2(1-\sigma^{\ast})-\frac{k-1}{2\a}}{2\a(3-2\sigma^{\ast}-\frac{1}{a_k})-1}.
\end{equation}
\end{thm}

It is an important problem to determine the exact value of $\sigma^{\ast}$. Generally it is a very difficult problem,
but if we assume the Lindel\"of hypothesis we can see easily that $\sigma^{\ast}=1 - 1/2a_k$. Hence from Theorem 1 we have

\begin{cor}
Suppose that the Lindel\"of hypothesis is true. If $(k-2)a_k < a_1+\cdots+a_{k-1}$, we have
\begin{equation*}
\int_1^T \Delta^2(\ba;x)dx=c(\ba)T^{1+\frac{k-1}{2\a}}+O\left(T^{1+\frac{k-1}{2\a}-\frac{2\a-(k-1)a_k}{2\a(2\a-1)a_k}+\varepsilon}\right).
\end{equation*}
where $c(\ba)$ is a certain positive constant.
\end{cor}

When $k=2$, we find  that $\sigma^{\ast}=1-1/2a_2$ unconditionally.
This is the consequence from the theorem on the  fourth power moment of the Riemann zeta function. Hence \eqref{maintheorem} gives

\begin{thm} In the case $k=2$, let $a_1 \leq a_2$. Then we have
\begin{equation} \label{cor1}
\int_1^T \Delta^2(a_1,a_2;x)dx=c_2 T^{1+\frac{1}{a_1+a_2}}+O\left(T^{1+\frac{1}{a_1+a_2}-\frac{a_1}{a_2(a_1+a_2)(a_1+a_2-1)}+\varepsilon}\right),
\end{equation}
where $c_2$ is a certain positive constant.
\end{thm}

This theorem improves the error term of  \eqref{ZhaiCao}.
We note that if we let $a_1=a_2=1$
the error term in \eqref{cor1} is $O(T^{1+\varepsilon})$, so  it can be said that
\eqref{cor1} corresponds to the result of Tong \eqref{Tong} except $T^{\varepsilon}$ factor.

Another interesting case is $k=3$.  We can show that
\begin{thm}
Let $k=3$. If $a_1 \leq a_2 \leq a_3$ and $a_3<a_1+a_2$, then we have
\begin{align*}
\int_1^T \Delta^2(a_1,a_2,a_3 ;x) dx=c_3T^{1+\frac{2}{a_1+a_2+a_3}}+O(T^{1+\frac{2}{a_1+a_2+a_3}-\eta_3+\varepsilon}),
\end{align*}
where
\begin{equation*}
\eta_3 =\begin{cases} \d \frac{1}{(a_1+a_2+a_3)(3+2(a_1+a_2+a_3)(1-1/a_3))} \qquad \mbox{if \ $3(a_2+a_3) \leq 7a_1$,} \\[1em]
                      \d \frac{4a_1a_3}{(a_1+a_2+a_3)\bigl((a_1+a_2+a_3)(a_1+3a_2+3a_3)(a_3-1)+a_3(5a_1+3a_2+3a_3)\bigr)}  \\[1em]
                        \hspace{3cm} \mbox{if \ $3(a_2+a_3)>7a_1, 3a_3+a_1 \leq 5a_2$ and $3a_3<a_1+3a_2$,}  \\[1em]
                      \d \frac{a_1+a_2-a_3}{a_3(a_1+a_2+a_3)(a_1+a_2+a_3-1)} \qquad \mbox{otherwise,}
       \end{cases}
\end{equation*}
and $c_3$ is a certain positive constant.
\end{thm}

We prove Theorem 3 in Section 4.


\section{The Tong-type identity of $\Delta(\ba;x)$ }


In \cite{Tong}, Tong studied the mean square of $\Delta(\underbrace{1,\ldots, 1}_{k};x)$. By using the functional equation of $\zeta^k(s)$
he derived a very useful identity of $\Delta(1,\ldots,1;x)$, which we call the Tong-type identity,
where the first finite sum is the same as that of the truncated Vorono\"i formula, while its error term is represented by the integrals
like \eqref{tongintegral} below.

In our case, by the functional equation of the Riemann zeta function
$$
\pi^{-s/2}\Gamma\left(\frac{s}{2}\right)\zeta(s)=\pi^{-(1-s)/2}\Gamma\left(\frac{1-s}{2}\right)\zeta(1-s),
$$
we find easily that the functional equation of $\varphi(s)$ and $\psi(s)$ has a form
\begin{align}
\Delta_1(s)\varphi(s)=\Delta_2(1-s)\psi(1-s),  \label{FE}
\end{align}
where
\begin{equation}
\Delta_1(s):=\prod_{j=1}^{k}\Gamma\left(\frac{a_j s}{2}\right)  \label{delta-1}
\end{equation}
and
\begin{equation}
\Delta_2(s):=\prod_{j=1}^{k}\Gamma\left(\frac{a_js-a_j+1}{2}\right).  \label{delta-2}
\end{equation}
Note that $\hat{d}(\ba;n)$ does not satisfy the Ramanujan conjecture and also the gamma factors on the left and right hand side of \eqref{FE}
are not the same for general $\ba$, so the pair of Dirichlet series $\varphi(s)$ and $\psi(s)$ is not contained in the so-called Selberg class.
In our previous paper \cite{CTZ}, we developed the theory of the Tong-type identity of the error term for such a pair of Dirichlet series.
In fact, we assumed that the two Dirchlet series $\varphi(s)$ and $\psi(s)$ satisfy the functional equation of the form
$\Delta_1(s)\varphi(s)=\Delta_2(r-s)\psi(r-s)$, where $\Delta_j(s)$ is the product of Gamma functions,  and derived the Tong-type identity
of the error term of the $\varrho$-th Riesz mean of the coefficients of $\varphi(s)$.

In order to write up the Tong-type identity for $\Delta(\ba;x)$ in the present case,
we follow the same notation of \cite{CTZ}. From \eqref{delta-1} and \eqref{delta-2}, we have  (we define $\a$ again for its importance.)
\begin{align*}
\a &= \frac{a_1+\cdots+a_k}{2}, \ \ \ r=1,\\
\mu&=\frac{1-k}{2}, \ \ \ \mu'=\sum_j \left(-\frac{a_j}{2}\right)+\frac12=-\a+\frac12, \\
\nu&=-\frac12\sum_j \log a_j, \ \ \ \nu'=-\frac12\sum_j a_j \log a_j, \\
\lambda&=\sum_j a_j \log a_j=\lambda',\\
h&=2\a e^{-\frac{\lambda+\lambda'}{2\a}}=2\a\prod_{j=1}^k a_j^{-a_j/\a}
\end{align*}
and
$$
\theta_{\varrho}=\frac{r}{2}-\frac{1}{4\a}+\varrho\left(1-\frac{1}{2\a}\right)+\frac{\mu'-\mu}{2\a}.
$$
In this paper we only consider the case $\varrho=0$, hence
\begin{align}
\theta_{0}&=\frac12-\frac{1}{4\a}+\frac{\mu'-\mu}{2\a}=\frac{k-1}{4\a} \label{theta0}.
\end{align}
We also need
\begin{align}
\lambda_0&=\theta_0+\frac{1}{2\a}-r-1=\frac{k+1}{4\a}-2.  \label{lambda0}
\end{align}

In Tong's theory, it is important to approximate $\Delta(\ba;x)$ by the $K$-th averaging integral
$$
\int_{\mathbf{E}_K}\Delta(\ba;\tilde{y})d\mathrm{Y}_K,
$$
where we use the notation 
\begin{equation*}
\int_{\mathbf{E}_K}g(\tilde{y})d\mathrm{Y}_K=\int_0^1\cdots \int_0^1 g(\tilde{y})dy_1\cdots dy_K,
\end{equation*}
with
$$
\tilde{y}=y+\frac{1}{x}(y_1+\cdots+y_K)
$$
for an integrable function $g(y)$.  
Let $\hat{\Delta}(\ba;x)$ be the
error term of the summatory function of $\hat{d}(\ba;n)$, which is defined by the same way as $\Delta(\ba;x)$.
Then the averaging integral can be expressed by the function defined by 
\begin{equation}  \label{tongintegral}
I(\lambda,M,N,y)=2\pi i \int_M^Nu^\lambda \hat{\Delta}(\ba;u)\exp\left(-ih (uy)^{\frac{1}{2\a}}\right)du.
\end{equation}

The next lemma gives the Tong-type identity of  $\Delta(\ba;y)$.
\begin{lem}  \label{Tong-identity}
Let $1 \leq x \leq y \leq (1+\delta)x$, $N=[x^{4\a-1-\varepsilon}]$ and $J=[(4\a^2r+4\a)\varepsilon^{-1}]$,
where $\delta$ is a small positive constant.
In every subinterval $[t,t+Bt^{1-1/2\a}] \subset [1,\sqrt{N}]$, there exists $M \neq \mu_n$ such that the following Tong-type
identity holds:
\begin{equation*}
\Delta(\ba;y)=\sum_{j=1}^7 R_{j}(y),
\end{equation*}
where
\begin{align*}
R_1(y)&=\kappa_0 y^{\theta_0}\sum_{\mu_n \leq M}\frac{b(n)}{\mu_n^{1-\theta_0}}\cos(h(y\mu_n)^{1/2\a}+c_0 \pi)   \\
   &=\kappa_0 \pi^{2\a(\theta_0-1)}y^{\theta_0}\sum_{n\leq M'}\frac{b(n)}{n^{1-\theta_0}}\cos(h\pi(yn)^{1/2\a}+c_0 \pi)   \\
   &=\kappa_0 \pi^{2\a\theta_0-k/2}y^{\theta_0}\sum_{n\leq M'}\frac{\hat{d}(\ba;n)}{n^{1-\theta_0}}\cos(h\pi(yn)^{1/2\a}+c_0 \pi),   \\
R_2(y)&=y^{\theta_0+\frac{1}{2\a}}\Re\{c_{00} I(\lambda_0, M, N,y)\}, \\
R_3(y)&=\Jsum_{l+m >0} \Re\left\{c_{lm} I\left(\lambda_0+\frac{l-m}{2\a}, M, N, y\right)\right\}x^{-l}y^{-l+\theta_0+\frac{1}{2\a}+\frac{l-m}{2\a}}, \\
R_4(y)&=\sum_{j=0}^K \sum_{m=0}^K \Re\left\{c_{jm}' I\left(\lambda_0-\frac{K+m}{2\a}, N, \infty, y+\frac{j}{x}\right)\right\}
        x^{K}\left(y+\frac{j}{x}\right)^{K+\theta_0+\frac{1}{2\a}-\frac{K+m}{2\a}}, \\
R_5(y)&= x^{\frac{k-3}{4\a}}M^{\max(\frac{k-3}{4\a},0)+\varepsilon}+x^{\frac{k+1}{4\a-2}}M^{\frac{k+1}{4\a}+\varepsilon}
         +x^{\frac{k-1}{4\a}-\frac12}M^{\omega_1-\frac32+\frac{k-1}{4\a}}\\
      & \quad +x^{(4\a-1)(1+\omega_1)-2K+\frac{k}{2\a}+\frac{2K}{\a}-6\a}, \\
R_6(y)&= 0, \\ 
R_7(y)&= \Delta(\ba;y)-\int_{{\mathbf{E}_K}}\Delta(\ba;\tilde{y})d\mathrm{Y}_K,
\end{align*}
where $M'=M/\pi^{2\a}$ and $\kappa_0\neq 0, c_{00}, c_{lm}, c_{jm}'$ are certain constants, $K$ is a suitably large integer and $\omega_1<1$ is
a certain constant.
\end{lem}

 This is Theorem 7 of \cite{CTZ}. We need one remark on $R_6(y)$. In fact $R_6(y)$ is given as
\begin{align*}
R_6(y) \ll \begin{cases} 0  & \mbox{if $b(n) \geq 0$}, \\
                         x^{\theta_0}M^{\omega_0-1+\frac{k-1}{4\a}}, & \mbox{if $b(n) \ll n^{\omega_0}$}.
            \end{cases}
\end{align*}
In our case we can take $R_6(y)=0$ since $b(n)=\pi^{2\a-k/2}\hat{d}(\ba,n)$ is always positive.

\medskip

We recall important evaluations of the integral of $I(\lambda,M,N,y)$ which we need in the next section.

\begin{lem} \label{I-sekibun}
Let $M <N<x^A$ ($A$ is a fixed positive number), $w$ be a real number and $0<\mu<\frac{M}{2}$. Then we have
\begin{align*}
&\int_x^{(1+\delta)x}I(\lambda,M,N,y) y^w \cos(h(\mu y)^{1/2\a}+c_0\pi)dy \\
& \hspace{2cm} \ll x^{w+1-3/4\a+\varepsilon} \max_{M \leq P \leq N}P^{\lambda+\sigma^{\ast}+1-3/4\a}.
\end{align*}
\end{lem}

\medskip

\begin{lem}  \label{Imeansquare-1}
Let $2(\lambda+\sigma^{\ast})\neq -1, \ M<N<x^A$ ($A$ is a fixed positive number) and $\delta>0$ with
$(1+\delta)^{1/\a}-1<1/4$. Then we have
$$
\int_x^{(1+\delta)x} |I(\lambda, M, N, y)|^2dy \ll x^{1-1/\a+\varepsilon} \max_{M \leq P \leq N}P^{2(\lambda+\sigma^{\ast}+1)-1/\a}.
$$
\end{lem}

\medskip

\begin{lem}  \label{Imeansquare-2}
Let $2(\lambda+\sigma^{\ast})\neq -1,  2(\lambda+\sigma^{\ast}+1)<1/\a,   M \geq 1$  and $\delta>0$ with
$(1+\delta)^{1/\a}-1<1/4$. Then we have
$$
\int_x^{(1+\delta)x} |I(\lambda, M, \infty, y)|^2dy \ll x^{1-1/\a+\varepsilon} M^{2(\lambda+\sigma^{\ast}+1)-1/\a}.
$$
\end{lem}

These lemmas are Lemmas 8, 9 and 10 of \cite{CTZ}, respectively.
See \cite{CTZ} for details.


\section{Mean square of $\Delta(\ba,x)$}


In the asymmetric many dimensional divisor problem, the number $(\mu'-\mu)/2=-\a+k/2$ plays an important role.
Though the proof of Theorem 1 goes in a similar way to \cite[Theorem 8]{CTZ},  we shall give a detailed proof for the sake of completeness.

Let
$$
K_1(y)=R_1(y)+R_2(y)
$$
and
$$
K_2(y)=\sum_{j=3}^7 R_j(y).
$$
It is sufficient to evaluate the integral $\int_x^{(1+\delta)x}(K_1(y)+K_2(y))^2dy$ for $1 \leq x<T$,
where $\delta$ is some fixed positive number.

We need the upper bound of the sum of $\hat{d} \,^2(\ba,n)$. Moreover we can prove
\begin{lem}   \label{d-hat-meansquare}
Let $x >1$. Then we have
\begin{equation}  \label{dhatsquare}
x^{2-1/a_k} \ll \sum_{n \leq x} \hat{d}\,^2(\ba;n) \ll x^{2-1/a_k+\varepsilon}.
\end{equation}
\end{lem}

\proof
By Cauchy's inequality we get
\begin{align*}
\hat{d}\,^2(\ba;n)&=\left(\sum_{n_1^{a_1}\cdots n_k^{a_k}=n}n_1^{a_1-1}\cdots n_k^{a_k-1}\right)^2 \\
&\leq \sum_{n_1^{a_1}\cdots n_k^{a_k}=n}1 \times \sum_{n_1^{a_1}\cdots n_k^{a_k}=n} n_1^{2(a_1-1)}\cdots n_k^{2(a_k-1)}\\
&\ll n^{\varepsilon}c(\ba;n),
\end{align*}
where $c(\ba;n)=\sum_{n_1^{a_1}\cdots n_k^{a_k}=n} n_1^{2(a_1-1)}\cdots n_k^{2(a_k-1)}$. We also note that $\hat{d}\,^2(\ba;n) \geq c(\ba;n)$.
The generating Dirichlet series of $c(\ba;n)$ has a form
$$
\sum_{n=1}^{\infty}\frac{c(\ba;n)}{n^s}=\prod_{j=1}^{k}\zeta(a_js-2(a_j-1)),  \quad \Re(s)>2-1/a_k.
$$
This Dirichlet series has poles at points $2-1/a_j \ \ (j=1, \ldots, k)$, hence
$$
\sum_{n \leq x}c(\ba;n) =c x^{2-1/a_k} \log^{A-1} x \cdot(1+o(1))
$$
where $c$ is some constant and $A$ is the number of $j$ such that $a_j=a_k$, hence the result follows.
\qed

\bigskip

Let $\sigma^{\ast}$ be the number defined by \eqref{def-sigmastar}. Assume that $\sigma^{\ast}$ satisfies \eqref{sigmastar-katei}.
The inequality \eqref{sigmastar-katei} is equivalent to
\begin{equation}
2(\lambda_0+\sigma^{\ast}+1)<\frac{1}{\a}, \label{katei-betsuhyouji}
\end{equation}
where $\lambda_0$ is defined by \eqref{lambda0}.


\subsection{Evaluation of $\int_x^{(1+\delta)x}K_1^2(y)dy$}


Let $\kappa_0'=\kappa_0 \pi^{2\a(\theta_0-1)}$ for simplicity.  By using
$$
\cos(x)\cos(y)=\frac12(\cos(x-y)+\cos(x+y))
$$
and substituting  \eqref{theta0} to $\theta_0$ in the Tong-type identity, we have
\begin{align*}
R_1(y)^2
&=\frac{\kappa_0'{}^2}{2}y^{\frac{k-1}{2\a}}\sum_{n \leq M'}\sum_{m \leq M'}\frac{b(n)b(m)}{(nm)^{1-\frac{k-1}{4\a}}}
     \left(\cos(h\pi y^{1/2\a}(n^{1/2\a}-m^{1/2\a}))  \right. \\
&  \hspace{4cm}       \left. +\cos(h\pi y^{1/2\a}(n^{1/2\a}+m^{1/2\a})+2c_0\pi)\right)\\
&=\frac{\kappa_0'{}^2}{2}\left(W_1(y)+W_2(y)+W_3(y)\right),
 \end{align*}
where
\begin{align*}
W_1(y)&=y^{\frac{k-1}{2\a}}\sum_{n \leq M'}\frac{b(n)^2}{n^{2-\frac{k-1}{2\a}}}, \\
W_2(y)&=y^{\frac{k-1}{2\a}}\Dsum_{\substack{n,m \leq M' \\ n \neq m}}\frac{b(n)b(m)}{(nm)^{1-\frac{k-1}{4\a}}}
                   \cos(h \pi y^{1/2\a}(n^{1/2\a}-m^{1/2\a})), \\
W_3(y)&=y^{\frac{k-1}{2\a}}\Dsum_{n,m \leq M' }\frac{b(n)b(m)}{(nm)^{1-\frac{k-1}{4\a}}}
                   \cos(h \pi y^{1/2\a}(n^{1/2\a}+m^{1/2\a})+2c_0\pi).
\end{align*}

For the integral of $W_1(y)$, we have
\begin{align*}
\int_x^{(1+\delta)x}W_1(y)dy = \sum_{n \leq M'}\frac{b(n)^2}{n^{2-\frac{k-1}{2\a}}} \int_x^{(1+\delta)x}y^{\frac{k-1}{2\a}}dy.
\end{align*}
Since \eqref{a-nokatei} is equivalent to $\frac{k-1}{2\a}<\frac{1}{a_k}$, we find that the series $\sum_{n}^{\infty}\frac{b(n)^2}{n^{2-\frac{k-1}{2\a}}}$
is convergent, hence by using \eqref{dhatsquare}, we have
\begin{equation*}
\sum_{n \leq M'}\frac{b(n)^2}{n^{2-\frac{k-1}{2\a}}}=\sum_{n=1}^{\infty} \frac{b(n)^2}{n^{2-\frac{k-1}{2\a}}}+O(M^{\frac{k-1}{2a}-\frac{1}{a_k}+\varepsilon}).
\end{equation*}
Hence
\begin{align} \label{W1-hyouka}
\int_x^{(1+\delta)x}W_1(y)dy = \sum_{n=1}^{\infty}\frac{b(n)^2}{n^{2-\frac{k-1}{2\a}}} \int_x^{(1+\delta)x}y^{\frac{k-1}{2\a}}dy
 +O(x^{1+\frac{k-1}{2\a}}M^{\frac{k-1}{2\a}-\frac{1}{a_k}+\varepsilon}).
\end{align}

By the first derivative test, we have
\begin{align*}
\int_x^{(1+\delta)x}W_2(y)dy & \ll x^{\frac{k-1}{2\a}+1-\frac{1}{2\a}} \Dsum_{\substack{m,n \leq M' \\ m \neq n}}\frac{b(n)b(m)}{(nm)^{1-\frac{k-1}{4\a}}}
\frac{1}{|n^{1/2\a}-m^{1/2\a}|} \\
&=x^{\frac{k-2}{2\a}+1}\left\{\Sigma_1+\Sigma_2\right\},
\end{align*}
where the summation conditions of $\Sigma_1$ and $\Sigma_2$ are given by
\begin{align*}
SC(\Sigma_1) : & \  |n^{1/2\a}-m^{1/2\a}| >\frac{1}{10}(nm)^{1/4\a} \\
\intertext{and}
SC(\Sigma_2) : & \  |n^{1/2\a}-m^{1/2\a}| \leq \frac{1}{10}(nm)^{1/4\a},
\end{align*}
respectively. It is easily seen that
\begin{align*}
\Sigma_1 &\ll \Dsum_{\substack{n,m \leq M' \\  |n^{1/2\a}-m^{1/2\a}| >\frac{1}{10}(nm)^{1/4\a} }}
              \frac{b(n)b(m)}{(nm)^{1-\frac{k-1}{4\a}}}\frac{1}{(nm)^{\frac{1}{4\a}}} \notag \\
         &\ll \left(\sum_{n \leq M'}\frac{b(n)}{n^{1-\frac{k-2}{4\a}}}\right)^2 \ll M^{\frac{k-2}{2\a}+\varepsilon},  
\end{align*}
where we have used the trivial estimate $\sum_{n \leq x}b(n) \ll x^{1+\varepsilon}$.
Next we consider $\Sigma_2$. By Lagrange's mean value theorem we have $n^{1/2\a}-m^{1/2\a}=\frac{1}{2\a}u_0^{1/2\a-1}(n-m)$ with some $u_0$ between
$n$ and $m$. Since $n \asymp m$ by $SC(\Sigma_2)$, we find
$$
|n^{1/2\a}-m^{1/2\a}| \geq (nm)^{1/4\a-1/2}|n-m|,
$$
thus we get
\begin{align*}
\Sigma_2 & \ll \Dsum_{\substack{n,m \leq M' \\ n \neq m}} \frac{b(n)b(m)}{(nm)^{\frac12-\frac{k-2}{4\a}}}\ \frac{1}{|n-m|} \\
& \ll \Dsum_{\substack{n,m \leq M' \\ n \neq m}} \left\{\left(\frac{b(n)}{n^{\frac12-\frac{k-2}{4\a}}}\right)^2+
       \left(\frac{b(m)}{m^{\frac12-\frac{k-2}{4\a}}}\right)^2\right\}\frac{1}{|n-m|}.
\end{align*}
By the symmetry on $n$ and $m$ and by using Lemma \ref{d-hat-meansquare} we obtain
\begin{align*}
\Sigma_2 \ll \Dsum_{\substack{n,m \leq M' \\ n \neq m}} \frac{b(n)^2}{n^{1-\frac{k-2}{2\a}}} \frac{1}{|n-m|}
         \ll M^{1-\frac{1}{a_k}+\frac{k-2}{2\a}+\varepsilon}.
\end{align*}
Here we note that the exponent of $M$ is $1-1/a_k+(k-2)/2\a \geq 0$ and $\Sigma_2$ is greater than $\Sigma_1$.
Hence
\begin{equation}  \label{W2-hyouka}
\int_x^{(1+\delta)x}W_2(y)dy \ll x^{\frac{k-2}{2\a}+1}M^{1-\frac{1}{a_k}+\frac{k-2}{2\a}+\varepsilon}.
\end{equation}
It is easily seen that $\int_x^{(1+\delta)x}W_3(y)dy$ is absorbed in the right hand side of \eqref{W2-hyouka}.

From \eqref{W1-hyouka} and \eqref{W2-hyouka}, we get
\begin{align}
\int_x^{(1+\delta)x}R_1^2(y)dy & = \frac{\kappa_0'{}^2}{2}\sum_{n=1}^{\infty}\frac{b(n)^2}{n^{2-\frac{k-1}{2\a}}}\int_x^{(1+\delta)x}y^{\frac{k-1}{2\a}}dy \notag \\
& \quad +O\left(x^{\frac{k-1}{2\a}+1+\varepsilon}M^{\frac{k-1}{2\a}-\frac{1}{a_k}}\right)
        +O\left(x^{\frac{k-2}{2\a}+1+\varepsilon}M^{\frac{k-2}{2\a}+1-\frac{1}{a_k}}\right). \label{R1}
\end{align}

Next we consider the mean square of $R_2(y)$. By Cauchy's inequality and Lemma \ref{Imeansquare-1}, we have
\begin{align*}
\int_x^{(1+\delta)x}R_2^2(y)dy & \ll x^{\frac{k-1}{2\a}+\frac{1}{\a}} \int_x^{(1+\delta)x}|I(\lambda_0,M,N,y)|^2dy \\
& \ll x^{\frac{k-1}{2\a}+\frac{1}{\a}} x^{1-\frac{1}{\a}+\varepsilon} \max_{M \leq P \leq N}P^{2(\lambda_0+\sigma^{\ast}+1)-\frac{1}{\a}}.
\end{align*}
We note that from \eqref{lambda0} and the assumption \eqref{sigmastar-katei}, $2(\lambda_0+\sigma^{\ast}+1)-1/\a <-1/a_k+(k-1)/2\a <0$. Therefore
\begin{equation}
\int_x^{(1+\delta)x}R_2^2(y)dy \ll x^{\frac{k-1}{2\a}+1+\varepsilon}M^{2\sigma^{\ast}-2+\frac{k-1}{2\a}}.   \label{R2}
\end{equation}

Finally we treat $\int_x^{(1+\delta)x}R_1(y)R_2(y)dy$. From the definitions of $R_1(y)$ and $R_2(y)$, we have
\begin{align*}
& \int_x^{(1+\delta)x} R_1(y)R_2(y)dy   \\
&=\Re \kappa_0' c_{00} \int_x^{(1+\delta)x} y^{\frac{k}{2\a}}I(\lambda_0,M,N,y)\sum_{n \leq M'}\frac{b(n)}{n^{1-\frac{k-1}{4\a}}}
  \cos(h\pi(ny)^{1/2\a}+c_0\pi) dy  \\
&=Re \kappa_0' c_{00} (I_1+I_2),
\end{align*}
where
\begin{equation*}
I_1=\int_x^{(1+\delta)x}y^{\frac{k}{2\a}}I(\lambda_0,M,N,y)\sum_{n\leq M'/2}\frac{b(n)}{n^{1-\frac{k-1}{4\a}}} \cos(h\pi(ny)^{1/2\a}+c_0\pi)dy
\end{equation*}
and
\begin{equation*}
I_2=\int_x^{(1+\delta)x}y^{\frac{k}{2\a}}I(\lambda_0,M,N,y)\sum_{M'/2 <n\leq M'}\frac{b(n)}{n^{1-\frac{k-1}{4\a}}} \cos(h\pi(ny)^{1/2\a}+c_0\pi)dy.
\end{equation*}
By Lemma \ref{I-sekibun} we have
\begin{align*}
I_1 \ll \sum_{n \leq M'}\frac{b(n)}{n^{1-\frac{k-1}{4\a}}}x^{\frac{k}{2\a}+1-\frac{3}{4\a}+\varepsilon}
\max_{M \leq P \leq N} P^{\lambda_0+\sigma^{\ast}+1-\frac{3}{4\a}}.
\end{align*}
By the assumption \eqref{sigmastar-katei}, the exponent of $P$ is negative, hence by using $\sum_{n \leq x}b(n) \ll x^{1+\varepsilon}$, we get
\begin{align}
I_1 & \ll x^{\frac{2k-3}{4\a}+1+\varepsilon}M^{\lambda_0+\sigma^{\ast}+1-3/4\a}\sum_{n \leq M'/2}\frac{b(n)}{n^{1-\frac{k-1}{4\a}}} \notag \\
& \ll x^{\frac{2k-3}{4\a}+1+\varepsilon}M^{\sigma^{\ast}-1+\frac{2k-3}{4\a}}. \label{R1R2-1}
\end{align}
Note that we avoid using Cauchy's inequality for $I_1$. But we apply Cauchy's inequality for $I_2$. In fact, we have
\begin{equation} \label{I-2}
I_2 \ll x^{\frac{k}{2\a}}(V_1 V_2)^{1/2},
\end{equation}
where we put
\begin{equation*}
V_1 =\int_x^{(1+\delta)x}|I(\lambda_0, M.N,y)|^2dy
\end{equation*}
and
\begin{equation*}
V_2 =\int_x^{(1+\delta)x}\left| \sum_{M'/2 <n\leq M'}\frac{b(n)}{n^{1-\frac{k-1}{4\a}}} \cos(h\pi(ny)^{1/2\a}+c_0\pi) \right|^2 dy.
\end{equation*}
For $V_1$ we apply Lemma \ref{Imeansquare-1} and get
\begin{align}
V_1 
& \ll x^{1-\frac{1}{\a}+\varepsilon} M^{2\sigma^{\ast}-2+\frac{k-1}{2\a}}. \label{V1}
\end{align}
By the similar method of evaluating the mean square of $R_1(y)$ we get
\begin{equation}
V_2 \ll xM^{\frac{k-1}{2\a}-\frac{1}{a_k}+\varepsilon}+x^{1-\frac{1}{2\a}+\varepsilon}M^{1-\frac{1}{a_k}+\frac{k-2}{2\a}}.  \label{V2}
\end{equation}
By \eqref{I-2}, \eqref{V1} and \eqref{V2} we get
\begin{equation}
I_2 \ll x^{1+\frac{k-1}{2\a}+\varepsilon}M^{\sigma^{\ast}-1+\frac{k-1}{2\a}-\frac{1}{2a_k}}
       +x^{1+\frac{2k-3}{4\a}+\varepsilon}M^{\sigma^{\ast}-\frac12+\frac{2k-3}{4\a}-\frac{1}{2a_k}}.  \label{R1R2-2}
\end{equation}

The estimates \eqref{R1}, \eqref{R2}, \eqref{R1R2-1} and  \eqref{R1R2-2} give the asymptotic formula of $\int_x^{(1+\delta)x}K_1^2(y)dy$.
There are six error terms arising from these formulas.
By $1-1/2a_k \leq \sigma^{\ast}$ and
$$
x^{1+\frac{2k-3}{4\a}}M^{\sigma^{\ast}-\frac12+\frac{2k-3}{4\a}-\frac{1}{2a_k}}
=\left(x^{\frac{k-2}{2\a}+1}M^{\frac{k-2}{2\a}+1-\frac{1}{a_k}}\right)^{1/2}
 \left(x^{\frac{k-1}{2\a}+1}M^{2\sigma^{\ast}-2+\frac{k-1}{2\a}}\right)^{1/2},
$$
it is easy to see that the first error term of \eqref{R1}, \eqref{R1R2-1} and \eqref{R1R2-2} are contained in the second error term of \eqref{R1}
and \eqref{R2}, hence we get
\begin{align}
\int_x^{(1+\delta)x}K_1^2(y)dy
& = \frac{\kappa_0'{}^2}{2}\sum_{n=1}^{\infty}\frac{b(n)^2}{n^{2-\frac{k-1}{2\a}}}\int_x^{(1+\delta)x}y^{\frac{k-1}{2\a}}dy \notag \notag \\
&   +O\left(x^{\frac{k-2}{2\a}+1+\varepsilon}M^{\frac{k-2}{2\a}+1-\frac{1}{a_k}}\right)
    +O\left(x^{\frac{k-1}{2\a}+1+\varepsilon}M^{2\sigma^{\ast}-2+\frac{k-1}{2\a}}\right).  \label{K1}
\end{align}


\subsection{Evaluation of $\int_x^{(1+\delta)x}K_2^2(y)dy$}


We evaluate the upper bounds of $\int_x^{(1+\delta)x}R_j^2(y)dy$ $(j=3,\ldots,7)$.
By Cauchy's inequality and Lemma \ref{Imeansquare-1}, we have
\begin{align*}
\int_x^{(1+\delta)x}R_3^2(y)dy
& \ll \Dsum_{\substack{0 \leq l,m \leq J \\ l+m>0}}x^{-4l+\frac{k+1}{2\a}+\frac{l-m}{\a}}\int_x^{(1+\delta)x}\left|I(\lambda_0+\frac{l-m}{2\a},M,N,y)\right|^2dy \\
& \ll \Dsum_{\substack{0 \leq l,m \leq J \\ l+m>0}}x^{-4l+\frac{k+1}{2\a}+\frac{l-m}{\a}}x^{1-\frac{1}{\a}+\varepsilon}
      \max_{M \leq P \leq N}P^{2(\lambda_0+\frac{l-m}{2\a}+\sigma^{\ast}+1)-\frac{1}{\a}} \\
& = \Sigma_3+\Sigma_4,
\end{align*}
where the summation conditions are
$$
SC(\Sigma_3): 0 \leq l \leq m \leq J, \ l+m>0 \quad \mbox{and} \quad SC(\Sigma_4): 0 \leq m<l \leq J.
$$
Since we assume $2(\lambda_0+\sigma^{\ast}+1)<1/\a$, we find that
\begin{align*}
\Sigma_3 &\ll \sum_{\substack{0 \leq m \leq l \leq J \\ l+m>0}}x^{-4l+\frac{k-1}{2\a}+\frac{l-m}{\a}+1+\varepsilon}
              M^{2(\lambda_0+\sigma^{\ast}+1)-\frac{1}{\a}+\frac{l-m}{\a}} \\
&=x^{\frac{k-1}{2\a}+1+\varepsilon}M^{2(\sigma^{\ast}-1)+\frac{k-1}{2\a}}\sum_{\substack{0 \leq m \leq l \leq J \\ l+m>0}}
 x^{-4l+\frac{l-m}{\a}}M^{\frac{l-m}{\a}}.
\end{align*}
The sum over $l$ and $m$ on the right hand side are bounded by
$$
\ll (xM)^{-1/\a}+x^{-4} \ll (xM)^{-1/\a},
$$
hence we have
\begin{equation} \label{Sigma3}
\Sigma_3 \ll x^{\frac{k-3}{2\a}+1+\varepsilon}M^{2(\sigma^{\ast}-1)+\frac{k-3}{2\a}}.
\end{equation}
Next we treat $\Sigma_4$. Since
$$
2(\lambda_0+\sigma^{\ast}+\frac{l-m}{2\a}+1)-\frac{1}{\a} \geq \frac{(a_k-a_1)+\cdots +(a_k-a_{k-1})+a_k}{(a_1+\cdots+a_k)a_k}>0
$$
in this case, we have
\begin{align*}
\Sigma_4 & \ll \sum_{0 \leq m < l \leq J}x^{-4l+\frac{k-1}{2\a}+1+\frac{l-m}{\a}+\varepsilon} N^{2(\lambda_0+\sigma^{\ast}+\frac{l-m}{2\a}+1)-\frac{1}{\a}} \\
&= x^{\frac{k-1}{2\a}+1+\varepsilon}N^{2(\lambda_0+\sigma^{\ast}+1)-\frac{1}{\a}}\sum_{0 \leq m < l \leq J}x^{-4l+\frac{l-m}{\a}}N^{\frac{l-m}{\a}}.
\end{align*}
Noting that $N=[x^{4\a-1-\varepsilon}]$, we find that the sum over $l$ and $m$ is $O(1)$, hence we have
\begin{equation}
\Sigma_4 \ll x^{\frac{k-1}{2\a}+1+\varepsilon}N^{2(\lambda_0+\sigma^{\ast}+1)-\frac{1}{\a}}. \label{Sigma4}
\end{equation}
From \eqref{Sigma3}, \eqref{Sigma4} and the assumption $M \leq \sqrt{N}$ we get
\begin{equation}  \label{R3}
\int_x^{(1+\delta)x}R_3^2(y)dy
\ll x^{\frac{k-3}{2\a}+1+\varepsilon}M^{2(\sigma^{\ast}-1)+\frac{k-3}{2\a}}+x^{\frac{k-1}{2\a}+1+\varepsilon}M^{4(\sigma^{\ast}-1)+\frac{k-1}{\a}}.
\end{equation}

Similarly we have
\begin{align*}
\int_x^{(1+\delta)x}R_4^2(y)dy
& \ll \sum_{j,m=0}^K x^{4K+\frac{k-1}{2\a}+\frac{1}{\a}-\frac{K+m}{\a}}
      \int_x^{(1+\delta)x}\left|I\left(\lambda_0-\frac{K+m}{2\a},N,\infty,y+\frac{j}{x}\right) \right|^2dy  \\
& \ll \sum_{j,m=0}^K x^{4K+\frac{k-1}{2\a}+\frac{1}{\a}-\frac{K+m}{\a}}x^{1-\frac{1}{\a}+\varepsilon}N^{2(\lambda_0-\frac{K+m}{2\a}+\sigma^{\ast}+1)-\frac{1}{\a}} \\
&  =  x^{4K+\frac{k-1}{2\a}+1-\frac{K}{\a}+\varepsilon}N^{2(\lambda_0+\sigma^{\ast}+1)-\frac{1}{\a}-\frac{K}{\a}}\sum_{j,m=0}^{K}(xN)^{-m/\a}.
\end{align*}
Here we have used Lemma \ref{Imeansquare-2}. Since the sum over $j$ and $m$ are bounded. we get, by the definition of $N$
\begin{align}
\int_x^{(1+\delta)x}R_4^2(y)dy
& \ll x^{4K+\frac{k-1}{2\a}+1-\frac{K}{\a}+\varepsilon}N^{2(\lambda_0+\sigma^{\ast}+1)-\frac{1}{\a}}x^{-(4\a-1-\varepsilon)\frac{K}{\a}} \notag \\
& \ll x^{\frac{k-1}{2\a}+1-\frac{K}{\a}+\varepsilon}N^{2(\lambda_0+\sigma^{\ast}+1)-\frac{1}{\a}}. \label{R4}
 \end{align}

Now consider $R_5(y)$. By taking $K$ large, we have
$$
R_5(y) \ll x^{\frac{k-3}{2\a}}M^{\max(\frac{k-3}{4\a}, 0)+\varepsilon}+x^{\frac{k+1}{4\a}-2}M^{\frac{k+1}{4\a}+\varepsilon}
+x^{\frac{k-1}{4\a}-\frac12}M^{-\frac12+\frac{k-1}{4\a}}.
$$
It is easy to see that
$$
R_5(y) \ll \begin{cases} x^{-1/4\a} & \mbox{if $k=2$} \\
                         x^{\frac{k-3}{4\a}}M^{\frac{k-3}{4\a}} & \mbox{if $k \geq 3$ and $M \ll x^{2\a-1}$}.
           \end{cases}
$$
Hence
\begin{equation}  \label{R5}
\int_x^{(1+\delta)x}R_5^2(y)dy \ll \begin{cases} x^{1-1/2\a} & \mbox{if $k=2$} \\
                                                 x^{1+\frac{k-3}{2\a}}M^{\frac{k-3}{2\a}} & \mbox{if $k \geq 3$ and $M \ll x^{2\a-1}$}.
           \end{cases}
\end{equation}

By the choice of $M$, $R_6(y)=0$, so we don't need to consider its mean square.

By the same method as that of \cite{CTZ}, we have
\begin{equation} \label{R7}
\int_x^{(1+\delta)x}R_7^2(y)dy \ll x^{\varepsilon}.
\end{equation}

The first error term on the right hand side of \eqref{R3} is clearly contained in the right hand side of \eqref{R5}.
Hence from \eqref{R3}, \eqref{R4}, \eqref{R5} and \eqref{R7}
\begin{align}
\int_x^{(1+\delta)x}K_2^2(y)dy
&\ll 
x^{\frac{k-1}{2\a}+1+\varepsilon}M^{4(\sigma^{\ast}-1)+\frac{k-1}{\a}}  \notag \\
& \quad +\begin{cases} x^{1-1/2\a} & \mbox{if $k=2$} \\
                                                 x^{1+\frac{k-3}{2\a}}M^{\frac{k-3}{2\a}} & \mbox{if $k \geq 3$ and $M \ll x^{2\a-1}$}.
         \end{cases}  \label{K2meansquare}
\end{align}


\subsection{Proof of Theorem 1}


Choose $M$ so that the two error terms of $\int_x^{(1+\delta)x}K_1^2(y)dy$ have the same order (see \eqref{K1}), namely
\begin{equation} \label{choice}
x^{\frac{k-2}{2\a}+1}M^{\frac{k-2}{2\a}+1-\frac{1}{a_k}} \asymp x^{\frac{k-1}{2\a}+1}M^{2(\sigma^{\ast}-1)+\frac{k-1}{2\a}},
\end{equation}
hence
\begin{equation} \label{choiceM}
M \asymp x^{\frac{1}{2\a(3-2\sigma^{\ast}-1/a_k)-1}}.
\end{equation}
Clearly $M$ satisfies $M \ll x^{2\a-1} \ll \sqrt{N}.$
Therefore we have
\begin{align}
\int_x^{(1+\delta)x}K_1^2(y)dy &= \frac{\kappa_0'{}^2}{2}\sum_{n=1}^{\infty}\frac{b(n)^2}{n^{2-\frac{k-1}{2\a}}}\int_x^{(1+\delta)x}y^{\frac{k-1}{2\a}}dy
+O\left(x^{1+\frac{k-1}{2\a}-\eta(\ba)+\varepsilon}\right), \label{K1meansquare}
\end{align}
where $\eta(\ba)$ is given by \eqref{eta}.

 First we shall see that $\int_x^{(1+\delta)x}K_2^2(y)dy$ is contained in the error term on the right hand side of \eqref{K1meansquare}.
Consider the first error term in \eqref{K2meansquare}. Since the exponent of $M$ is negative, it is smaller than the error term of \eqref{K1meansquare}.
Next consider the second error term of \eqref{K2meansquare}. For $k=2$ there is nothing to prove. For $k \geq 3$, it is enough to show that
$$
\frac{1}{\a}-\frac{k-3}{2\a} \cdot \frac{1}{2\a(3-2\sigma^{\ast}-1/a_k)-1} > \eta(\ba),
$$
or equivalently $2-1/a_k>\sigma^{\ast}$. This is true from the assumption \eqref{sigmastar-katei}.

Finally we treat $\int_x^{(1+\delta)x}K_1(y)K_2(y)dy$. By Cauchy's inequality, \eqref{K1meansquare} and \eqref{K2meansquare}  we have
\begin{align*}
\int_x^{(1+\delta)x}K_1(y)K_2(y)dy
&\ll \left(\int_x^{(1+\delta)x}K_1^2(y)dy\right)^{1/2}\left(\int_x^{(1+\delta)x}K_2^2(y)dy\right)^{1/2} \\
& \ll x^{1+\frac{k-1}{2\a}+\varepsilon}M^{2(\sigma^{\ast}-1)+\frac{k-1}{2\a}}
      + \begin{cases} x & \mbox{if $k=2$} \\
                    x^{1+\frac{k-2}{2\a}}M^{\frac{k-3}{4\a}}  & \mbox{if $k \geq 3$}
        \end{cases}.
\end{align*}
Since $M$ is chosen by the equation \eqref{choice}, this is also contained in the error term of \eqref{K1meansquare}.
This completes the proof of Theorem 1.  \qed


\section{Proof of Theorem 3}


In order to prove Theorem 3 we need some preparation. Define $m(\sigma)$ ($1/2 \leq \sigma <1$) as the supremum of all numbers $m$ such that
$$
\int_1^T|\zeta(\sigma+it)|^m dt \ll T^{1+\varepsilon}.
$$
It is known that $ m(\sigma)\geq 4$ for $\sigma \geq 1/2$, $m(7/12) \geq 6$ and $m(5/8) \geq 8$.
Ivi\'c studied $m(\sigma)$ in great detail. One can find a lower bound in \cite[Theorem 8.4]{I1}.
Here we can suppose that $m(\sigma)$ is continuous.
Especially we have the following simpler but a little weaker form:
\begin{equation}  \label{m-sigma}
m(\sigma) \geq \begin{cases} \d \frac{4}{3-4\sigma} & \mbox{if \  $ \d \frac12 \leq \sigma \leq \frac58$} \\[1em]
                              \d \frac{3}{1-\sigma}  & \mbox{if \  $ \d \frac58 \leq \sigma <1$}.
                 \end{cases}
\end{equation}

The following lemma is used essentially in Ivi\'c's argument \cite{I2}.
\begin{lem}
Let $a_j \ (1 \leq j \leq k)$ be positive integers such that $a_1 \leq \cdots \leq a_k$ and let $\psi(s)$ and $\sigma^{\ast}$ be
defined by \eqref{def-psi} and \eqref{def-sigmastar}, respectively.  Define the function $H(\sigma)$ by
$$
H(\sigma)=\sum_{j=1}^k \frac{1}{m(a_j \sigma-a_j+1)}.
$$
If
\begin{equation*}
H(\sigma) \leq 1/2   
\end{equation*}
for some $\sigma$, we have $\sigma^{\ast} \leq \sigma$.
\end{lem}

\proof
We write $\sigma_j=a_j\sigma-a_j+1$ for simplicity. Suppose that
$$
\sum_{j=1}^k \frac{1}{m(\sigma_j)} \leq \frac12.
$$
Then by H\"older's inequality, we have
\begin{align*}
\int_1^T|\psi(s)|^2dt &=\int_1^T\prod_{j=1}^k |\zeta(\sigma_j+ia_jt)|^2dt \\
& \leq \prod_{j=1}^k \left(\int_1^T |\zeta(\sigma_j+ia_jt)|^{m(\sigma_j)}dt\right)^{\frac{2}{m(\sigma_j)}}
\left(\int_1^T 1 dt\right)^{1-\sum_{j=1}^k \frac{2}{m(\sigma_j)}}\\[1ex]
& \ll T^{1+\varepsilon}.
\end{align*}
Hence from the definition of $\sigma^{\ast}$, we have $\sigma^{\ast} \leq \sigma$.  \qed

\medskip

We remark that since $H(\sigma)$ is decreasing,  if we have
$$
H\left(1-\frac{k-1}{4\a}\right) < \frac12,
$$
then Theorem 1 holds.

\bigskip

\begin{lem}  Let $k=3$,  $a_1 \leq a_2 \leq a_3$ and $a_3<a_1+a_2$.  Let $\sigma^{\ast}$ be defined by
\eqref{def-sigmastar}. Then we have
\begin{align} \label{three-dim}
 \sigma^*
\begin{cases}
\le  \d 1-\frac {5}{4(a_1+a_2+a_3)} & \mbox{if $3(a_2+a_3)\le 7a_1$, }\\[1em]
\le  \d 1-\frac {3}{a_1+3a_2+3a_3} & \mbox{if $3(a_2+a_3)> 7a_1$, $3 a_3+a_1\le 5 a_2$ and $3 a_3< a_1+3 a_2$,}\\[1em]
=    \d 1-\frac {1}{2a_3} & \mbox{otherwise.}
\end{cases}
\end{align}
\end{lem}

\proof
We shall find $1-\frac{1}{2a_3} \leq \sigma <1-\frac{1}{a_1+a_2+a_3}$ such that $H(\sigma) \leq 1/2$.
For simplicity we put $\sigma_j=a_j\sigma-a_j+1 \ (j=1,2,3)$ for $\sigma \in [\frac 12 ,1]$ as before.
It is easy to see that $\frac 12 \le \sigma_3\le\sigma_2\le \sigma_1<1$.

Now we use a weak version \eqref{m-sigma} of Theorem 8.4 in Ivi\'c.

\medskip

(Case 1) We first consider the case $3(a_2+a_3)\le 7a_1$. Put
$$
\sigma:=1-\frac {5}{4(a_1+a_2+a_3)}.
$$
Clearly $\sigma <1-1/(a_1+a_2+a_3)$.  Since $3a_3\le 7a_1-3a_2\le (2a_1+5a_2)-3 a_2=2(a_1+a_2)$, we have $\sigma\ge 1-\frac{1}{2a_3}$
and $\sigma_1\le \frac 58$. By \eqref{m-sigma} we have
\begin{align*}
H(\sigma)&=\sum_{j=1}^3\frac {1}{m(\sigma_j)} \le \frac{3-4\sigma_1}{4}+\frac{3-4\sigma_2}{4}+\frac{3-4\sigma_3}{4}
=\frac 12.
\end{align*}
Hence we get $\sigma^{\ast} \leq \sigma$.

\medskip

(Case 2)  When $3(a_2+a_3)> 7a_1$, $3 a_3+a_1\le 5 a_2$ and $3 a_3< a_1+3 a_2$, we put
$$
\sigma:= 1-\frac {3}{a_1+3a_2+3a_3}.
$$
It is clear that $\sigma<1-1/(a_1+a_2+a_3)$ and  $\sigma>1-\frac{1}{2a_3}$ by the last condition.
One can check that by the first two conditions that $\frac 58<\sigma_1<1$ and $\frac 12\le \sigma_3\le \sigma_2\le \frac 58$ . Hence
\begin{align*}
H(\sigma)&=\sum_{j=1}^3\frac {1}{m(\sigma_j)}\\
&\le\frac{1-\sigma_1}{3}+\frac{3-4\sigma_2}{4}+\frac{3-4\sigma_3}{4}
=\frac 12.
\end{align*}
Hence we get $\sigma^{\ast} \leq \sigma$.

\medskip

(Case 3) We consider the case $3(a_2+a_3)> 7a_1$, $3 a_3+a_1\le 5 a_2$ and $3 a_3\ge a_1+3 a_2$. In this case we put
$$
\sigma:= 1-\frac {1}{2a_3}.
$$
Note that this is the best possible choice. We easily check  by the last condition that
\begin{align*}
3a_3\ge a_1+3a_2\ge 4a_1,
\end{align*}
hence
\begin{align*}
\sigma_1=a_1\left(1-\frac {1}{2a_3}\right)-a_1+1=1-\frac {a_1}{2a_3}\ge
\frac 58.
\end{align*}

Now we consider two cases.

(i) If $3a_3\le 4a_2$, then $\sigma_2\le \frac 58$. We use the third condition to get
\begin{align*}
H(\sigma)&=\sum_{j=1}^3\frac {1}{m(\sigma_j)}\le\frac{1-\sigma_1}{3}+\frac{3-4\sigma_2}{4}+\frac 14\\
&=\frac{a_1+3a_2}{6a_3}\le\frac 12.
\end{align*}

(ii) If $3a_3> 4a_2$, then $\sigma_2 > \frac 58$.
By the third condition we have $3a_3\ge a_1+3a_2\ge 2(a_1+a_2)$. Hence
\begin{align*}
H(\sigma)&=\sum_{j=1}^3\frac {1}{m(\sigma_j)}\le\frac{1-\sigma_1}{3}+\frac{1-\sigma_2}{3}+\frac 14\\
&=\frac 14+\frac{a_1+a_2}{6a_3}\le\frac 14+\frac 16 \cdot \frac 32=\frac 12.
\end{align*}
Combining the two cases (i) and (ii), we have $\sigma^{\ast}=\sigma=1-1/(2a_3)$.

\medskip

(Case 4) Finally we consider the case $3(a_2+a_3)> 7a_1$, $3 a_3+a_1> 5 a_2$. In this case we put
$$
\sigma:= 1-\frac {1}{2a_3}.
$$

In this case we easily check  by the second condition that
\begin{align*}
&3a_3>5 a_2-a_1\ge 4a_2,\\
\intertext{hence}
&\sigma_2=a_2\left(1-\frac {1}{2a_3}\right)-a_2+1=1-\frac {a_2}{2a_3}>\frac 58.
\end{align*}
We have $1>\sigma_2>\frac 58, 1>\sigma_1>\frac 58$,   and $\sigma_3=\frac 12$.
Hence
\begin{align*}
H(\sigma)&=\sum_{j=1}^3\frac {1}{m(\sigma_j)}\le\frac{1-\sigma_1}{3}+\frac{1-\sigma_2}{3}+\frac{1}{4}\\
&\le \frac 18+\frac 18+\frac 14=\frac 12.
\end{align*}
Therefore we have $\sigma^{\ast}=\sigma=1-1/(2a_3).$
\qed

\bigskip

\noindent{\it Proof of Theorem 3.}
Now the proof of Theorem 3 is immediate by substituting each value on the right hand side of \eqref{three-dim}
to \eqref{eta}.

\bigskip

{\bf Remark.}  From Lemma 7 we have
\begin{align*}\sigma^*(3,4,5)=\frac{9}{10},\ \ \sigma^*(2,3,4)=\frac 78,\end{align*}\mbox
which are the best possible results .
By Theorem 8.4 of Ivi\'c\cite{I1} we also note the following slightly better results
\begin{align*}
\sigma^*(4,5,6)\le \frac{214}{233},\ \
\sigma^*(1,2,2)\le \frac{41761}{54522}=0.765948\cdots .
\end{align*}

\bigskip

\begin{flushleft}
Xiaodong Cao \\
Department of Mathematics and Physics, \\
Beijing Institute of Petro-Chemical Technology,\\
Beijing, 102617, P. R. China \\
e-mail: caoxiaodong@bipt.edu.cn

\bigskip

Yoshio Tanigawa\\
Graduate School of Mathematics,\\
Nagoya University, \\
Nagoya, 464-8602, Japan\\
e-mail: tanigawa@math.nagoya-u.ac.jp

\bigskip

Wenguang Zhai\\
Department of Mathematics, \\
China University of Mining and Technology, \\
Beijing 100083, P. R. China\\
e-mail: zhaiwg@hotmail.com
\end{flushleft}


\medskip

\begin{thebibliography}{99}

\bibitem{CTZ1} X. Cao, Y. Tanigawa and W. Zhai, On a conjecture of Chowla and Walum, Science China Mathematics {\bf 53} (2010), 2755--2771.

\bibitem{CTZ} X. Cao, Y. Tanigawa and W. Zhai, On the Tong-type identity and the mean square of the error term for an
extended Selberg class, submitted.

\bibitem{CK} K. Corr\'adi and I. K\'atai, A comment of K. S. Ganggadharan's paper entitled ``Two classical lattice point
problems", Magyar Tud. Akad. Mat. Fiz.Oszt. K\"ozl. {\bf 17} (1967), 89--97.

\bibitem{C} H. Cram\'er, \"Uber zwei S\"atze von Hern G. H. Hardy, Math. Z. {\bf 15} (1922), 201--210.

\bibitem{Hafner} J. L. Hafner, New omega theorems for two classical lattice point problems, Invent. Math. {\bf 63} (1981), 181-186.

\bibitem{Huxley} M. N. Huxley, Exponential sums and lattice points III, Proc. London Math. Soc. {\bf 87} (2003), 591--609.

\bibitem{I1} A. Ivi\'c, The Riemann Zeta-Function, John Wiley and Sons, 1985.

\bibitem{I2} A. Ivi\'c, The general divisor problem, J. Number Theory {\bf 27} (1987), 73--91.

\bibitem{IKKN}A. Ivi\'c,  E. Kr\"atzel, M. K\"uhleitner and W. G. Nowak, Lattice points in large regions and
related arithmetic functions: recent developments in a very
classic topic. (English summary) Elementare und analytische
Zahlentheorie, 89--128, Schr.Wiss. Ges. Johann Wolfgang Goethe
Univ.Frankfurt am Main, 20, Franz Steiner Verlag Stuttgart,
Stuttgart, 2006.


\bibitem{Kra} E. Kr\"{a}tzel, Lattice Points, Kluwer Academic Publishers, Dordrecht 1988.

\bibitem{LT} Y.-K. Lau and K.-M. Tsang, On the mean square formula of the error term in the Dirichlet divisor problem, Math. Proc.
Camb. Phil. Soc. {\bf 146} (2009), 277--287.

\bibitem{Tong} K.-C. Tong, On divisor problem III, Acta Math. Sinica {\bf 6} (1956), 515--541.

\bibitem{ZC} W. Zhai and X. Cao, On the mean square of the error term for the
asymmeric two-dimensional divisor problem (I), Monatsh. Math. {\bf 159} (2010), 185--209.
\end{thebibliography}
\end{document}